\documentclass[english]{article}
\usepackage[T1]{fontenc}
\usepackage[latin9]{inputenc}
\pdfoutput=1
\usepackage{float}
\usepackage{amsmath}
\usepackage{amssymb}
\usepackage{subcaption}
\usepackage[labelformat=parens,labelsep=quad,skip=3pt]{caption}
\usepackage{graphicx}
\usepackage{setspace}
\usepackage{lineno}
\usepackage{multirow}
\usepackage{hyperref}
\usepackage{soul}
\hypersetup{
    colorlinks=true,
    linkcolor=blue,
    filecolor=magenta,      
    urlcolor=cyan,
}
\makeatletter

\floatstyle{ruled}
\newfloat{algorithm}{tbp}{loa}
\providecommand{\algorithmname}{Algorithm}
\floatname{algorithm}{\protect\algorithmname}

\usepackage{beamerarticle,pgf}

\newcommand\makebeamertitle{\frame{\maketitle}}
\AtBeginDocument{
	\let\origtableofcontents=\tableofcontents
	\def\tableofcontents{\@ifnextchar[{\origtableofcontents}{\gobbletableofcontents}}
	\def\gobbletableofcontents#1{\origtableofcontents}
}

\makeatother

\usepackage{babel}
\begin{document}

\title{A computational geometry method for the inverse scattering problem}
\author{Maria L. Daza-Torres\footnotemark[2]\ \footnotemark[3]
\and
Juan Antonio Infante del R\'io\footnotemark[4]
\and
Marcos A. Capistr\'an\footnotemark[2]
\and 
J. Andr\'es Christen\footnotemark[2]
}

\renewcommand{\thefootnote}{\fnsymbol{footnote}}

\footnotetext[2]{Centro de Investigaci\'on en Matem\'aticas (CIMAT), Jalisco S/N, Valenciana, Guanajuato,
36023, M\'exico.
\textit{marcos@cimat.mx}}
\footnotetext[4]{Instituto de Matem\'atica Interdisciplinar y Departamento de An\'alisis Matem\'atico y Matem\'atica Aplicada,
Facultad de CC. Matem\'aticas, Universidad Complutense de Madrid, Plaza de Ciencias 3,
28040, Madrid, Spain.
\textit{infante@mat.ucm.es}}

\footnotetext[3]{Corresponding author}

\makebeamertitle
In this paper we demonstrate a computational method to solve the inverse
scattering problem for a star-shaped, smooth, penetrable obstacle
in 2D. Our method is based on classical ideas from computational geometry.
First, we approximate the support of a scatterer by a point cloud.
Secondly, we use the Bayesian paradigm to model the joint conditional
probability distribution of the non-convex hull of the point cloud
and the constant refractive index of the scatterer
given near field data. Of note, we use the non-convex hull of the
point cloud as spline control points to evaluate, on a finer mesh,
the volume potential arising in the integral equation formulation
of the direct problem. Finally, in order to sample the arising posterior
distribution, we propose a probability transition kernel that commutes
with affine transformations of space. Our findings indicate that our
method is reliable to retrieve the support and constant refractive 
index of the scatterer simultaneously. Indeed, our sampling method is robust 
to estimate a quantity of interest such as the area of the scatterer. We conclude 
pointing out a series of generalizations of our method.

\section{Introduction}\label{sec:intro}

In this paper we demonstrate a computational method to solve the inverse
scattering problem for a star-shaped, smooth, penetrable obstacle
in 2D. Broadly speaking, inverse problems in physical systems are
statistical inference problems restricted by an initial and/or boundary value
problem for a partial differential equation. Consequently, inverse
problems are difficult since the quantity to be inferred is defined
in a function space, while data is finite and noisy. Stuart \cite{Stuart2010}
studied conditions for the well-posedness of the Bayesian formulation
of inverse problems. Recently, Christen \cite{christen2017posterior} showed
that that the Bayesian formulation of inverse problems is well posed under
more general conditions. In either case, Bayesian inverse problems defined on function 
spaces have two main difficulties: the computationally efficient representation of the 
function to be estimated and the sampling of a posterior measure defined on an infinite-dimensional space. 
Although there are many theoretical approaches to overcome this issues 
\cite{bui2013computational,Girolami2011,Stuart2010},
the development of competitive computational methods for inverse
problems remains a very active research field.

We rely on classical computational geometry methods
to deal with the issues of function representation and sampling. First,
instead of posing the problem in an infinite dimensional setting,
we approximate the support of each scatterer using a point cloud \cite{edelsbrunner1983shape}.
The rationale of our approach is that the direct problem has a finite
rank property \cite{birman1977estimates}, i.e., the Frechet derivative
of the mapping from scattering obstacle to near field is a compact
operator, and its singular values have zero as a limit point \cite{colton2012inverse}.
Secondly, we use the Bayesian paradigm to model the joint conditional
probability of the non-convex hull of the point cloud and the scatterer constant
refractive index given measurements of the near field (in the remainder
of this paper, we shall refer to this conditional probability as the
posterior distribution). Of note, we use the non-convex hull of the
point cloud as spline control points to evaluate, on a finer mesh,
the volume potential arising in the integral equation formulation
of the direct problem. Finally, in order to sample the arising posterior distribution, we
construct a probability transition kernel based on moves that commute
with affine transformations of space. This property, first introduced
by Christen and Fox \cite{christen2010}, implies that the
performance of the Markov Chain Monte Carlo (MCMC) method is independent
of the aspect ratio of anisotropic distributions. 

The outline of the paper is as follows. In Section 2, we introduce
the direct scattering problem and propose a numerical solution. Section 3
is devoted to the representation of the shape of an obstacle based
on a point cloud. This representation uses the $\alpha$-shape algorithm
and allows us to parameterize the obstacle shape with a parameter
and a cloud of points. In Section 4, we present the Bayesian formulation
of the inverse problem. The solution of the inverse problem is the
posterior probability distribution. Thereby, an affine invariant MCMC
method to sample the posterior distribution is designed. Numerical
results and discussion are presented in Section 5. Finally, Section 6
summarizes our findings, pointing out limitations and possible generalizations
of our inference method.

\section{The direct scattering problem}\label{sec:dir_prob}

In this Section we describe the direct scattering problem and propose a method
for its numerical solution. We refer the reader to \cite{colton2012inverse}
for a more detailed discussion about the mathematical formulation
of the scattering problem for a penetrable obstacle.

Let us assume that the scatterer $D\subset\mathbb{R}^{2}$ is bounded
and penetrable. If the incident wave is a plane wave, the acoustic
scattering problem can be cast as the following problem for the Helmholtz
equation (see \cite{colton2012inverse}),
\begin{equation}
\begin{split}\Delta u(\mathbf{x})+k^{2}(b(\mathbf{x})+1)u(\mathbf{x}) & =0\quad\text{in}\quad\mathbb{R}^{2},\\
u(\mathbf{x}) & =u^{i}(\mathbf{x})+u^{s}(\mathbf{x}),\\
\lim_{r\rightarrow\infty}r\left(\frac{\partial u^{s}(\mathbf{x})}{\partial r}-iku^{s}(\mathbf{x})\right) & =0,
\end{split}
\label{eq:Helmholtz}
\end{equation}
where $u$ is the total field, $u^{i}(\mathbf{x})=\text{e}^{ik\mathbf{x}\cdot d}$
is the incident field, $u^{s}(\mathbf{x})$ is the scattered field,
$r=|\mathbf{x}|$, and $k$ is the wavenumber. The refractive index
is assumed to be constant on $D$,

\[
b(\mathbf{x})=\left\{ \begin{array}{lcc}
{\bf b}, &  \mathbf{x}\in D\\
\\
0, &  \mathbf{x}\in\mathbb{R}^{2}-D
\end{array}\right.
\]
for ${\bf b}\in\mathbb{R},$ ${\bf {b}\neq-1}$.\\

The Helmholtz equation \label{The-Helmholtz-equation} is equivalent
to the Lippmann-Schwinger equation 
\begin{equation}
u(x)=u^{i}-k^{2}\int_{D}\Phi(k|x-y|)b(y)u(y)dy,\label{eq:LS}
\end{equation}
 where
\begin{equation}
\Phi(r)=-\frac{i}{4}H_{0}^{(1)}(r)
\end{equation}
 and $H_{0}^{(1)}$ is the Hankel function of the first class of order
zero.

\subsection{Numerical solution of the direct scattering problem}\label{sec:num_method}

In order to obtain an efficient solution of the forward model, we
approximate the volume potential in equation \eqref{eq:LS} using a corrected trapezoidal rule. 
The correction coefficients are calculated in \cite{aguilar2002high}.
Taking the first correction coefficient $c_{1}$, we obtain a fourth
order approximation

\begin{equation}
u_{j,h}=u_{j,h}^{i}-h^{2}k^{2}\sum_{l\in\mathbb{Z}_{N}^{2}}\Phi_{j-l,h}b_{l,h}u_{l,h}\hspace{3mm}+O(h^{4}),\hspace{2mm}j\in\mathbb{Z}_{N}^{2}\label{eq:trapezoidal}
\end{equation}
where $\mathbb{Z}_{N}^{2}=\lbrace j\in\mathbb{Z}^{2}:0\leq j_{k}\leq N,k=1,2\rbrace$,

\[
\Phi_{j,h}=\left\{ \begin{array}{lc}
\Phi(k|j|h) & j\neq0\\
\\
\Phi_{j,h}=\beta_{1} & j=0,
\end{array}\right.
\]
$\beta_{1}=-\frac{i}{4}+\frac{1}{2\pi}\left(\ln\left(\frac{hk}{2}\right)+\gamma+c_{1}\right)$,
and $\gamma=0.5772156649015328606\ldots$ is Euler\textquoteright s
constant. \\

The discrete approximation of equation \eqref{eq:LS} can be written in matrix form 

\begin{equation}
(K_{D,b}+I)u_{h}=u_{h}^{i}\label{eq:matrix_form_T}.
\end{equation}

\begin{figure}[H]
\begin{centering}
\includegraphics[scale=0.25]{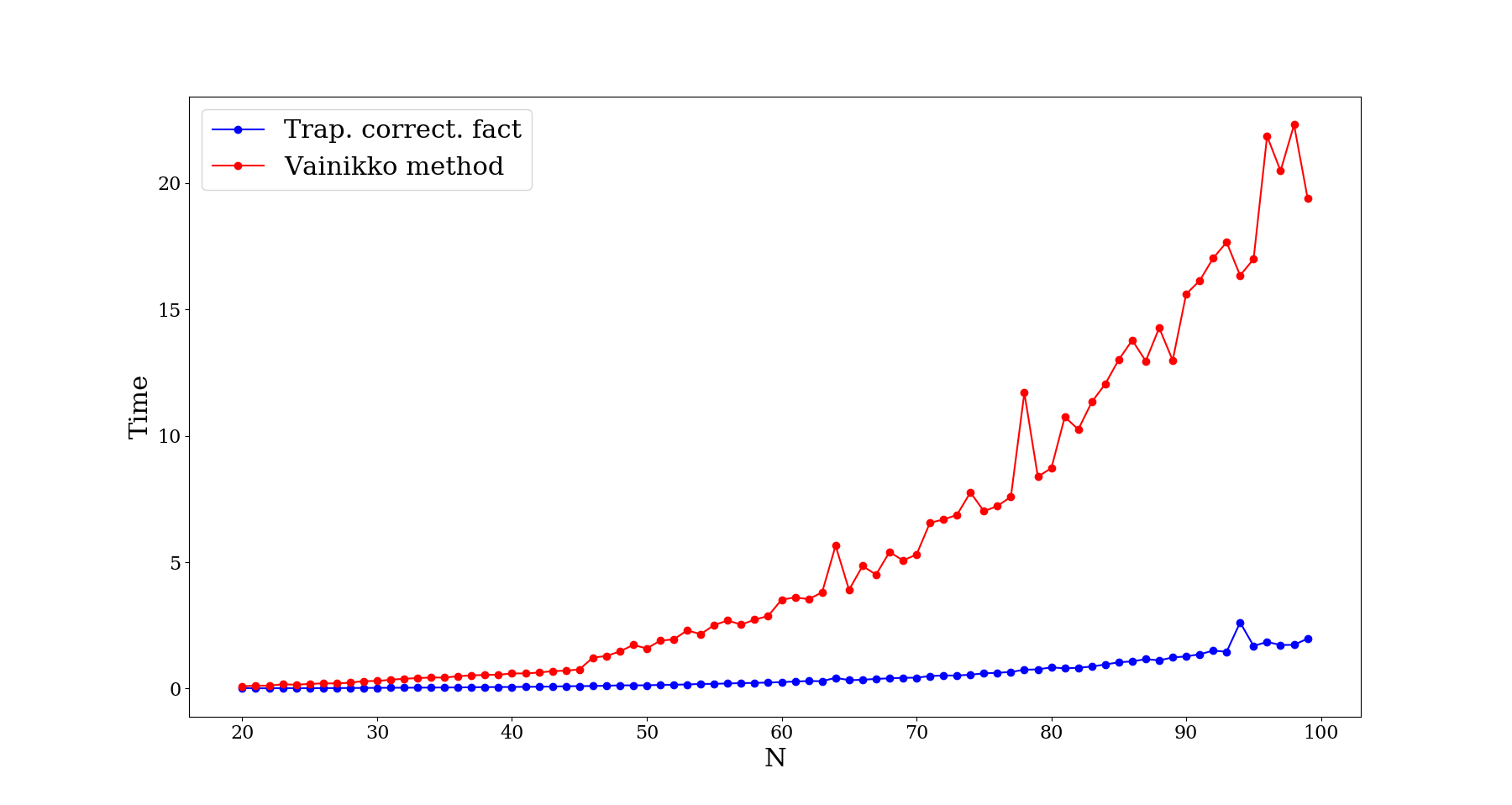}
\par\end{centering}
\caption{Comparison of computational times between (in seconds)  the method of Vainikko (red) 
and the method propose (blue) for solving the Lippmann-Schwinger equation 
with the configuration depicted in Example 1 of Section 5. \label{fig:gmresVSfact}}
\end{figure}

Note that $K_{D,b}$ is a sparse matrix that depends on the obstacle
and the refractive index. Of note, we can readily index the entries
where the matrix is different from zero using the definition of $b$.
Thus, the matrix structure allows us to solve the system efficiently 
(\ref{eq:matrix_form_T}). In Figure \ref{fig:gmresVSfact}, we compare
computational times to solve the Lippmann-Schwinger equation numerically
using equation \eqref{eq:matrix_form_T} with the numerical solution of Vainikko
\cite{vainikko2000fast}. It is apparent that the method introduced
above is considerably faster than the method of Vainikko.

\section{Point cloud approximation of a scatterer}\label{sec:point_cloud}

\begin{figure}[H]
\begin{centering}
\includegraphics[scale=0.25]{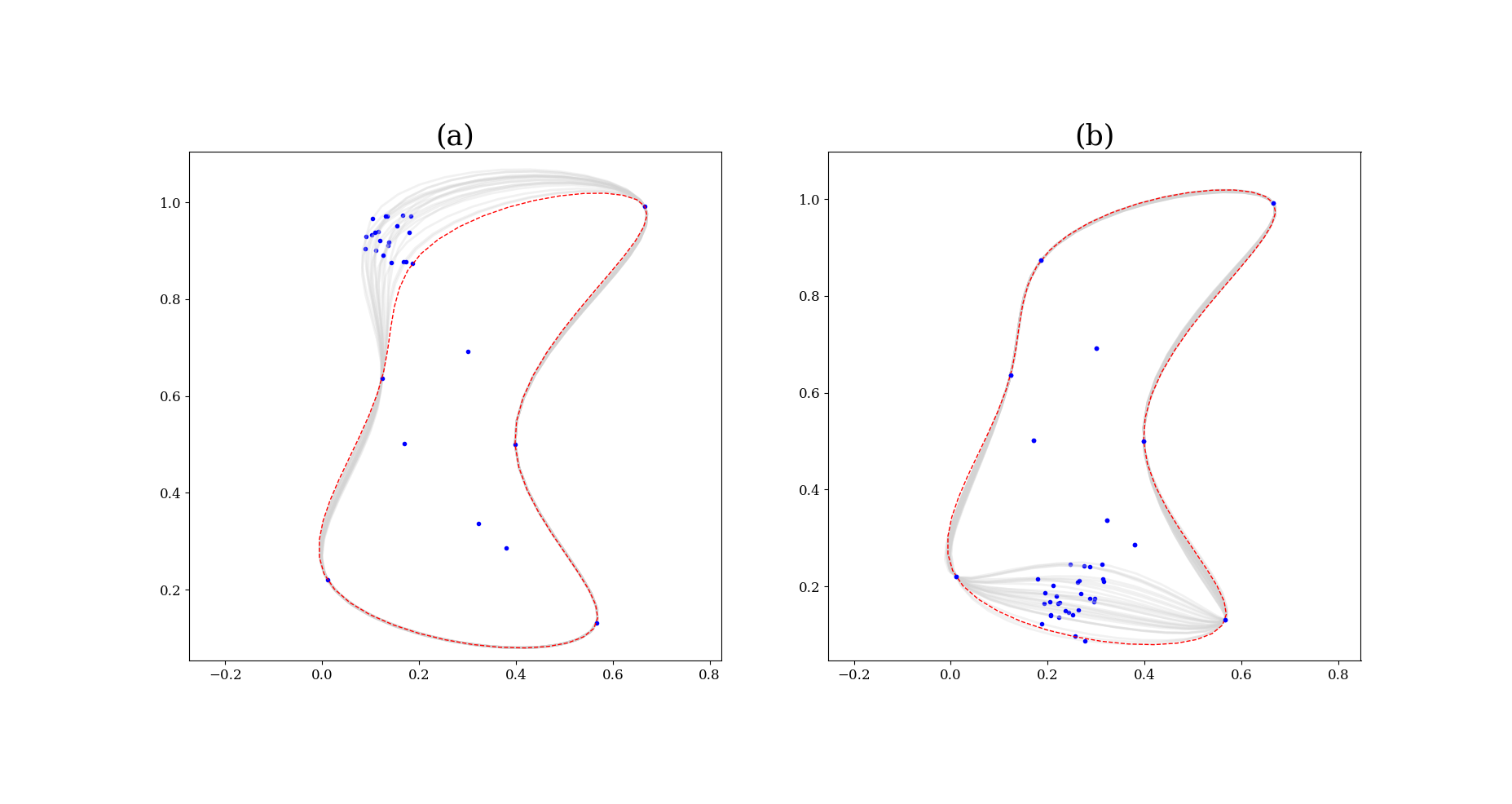}
\par\end{centering}
\caption{Part (a) shows the interpolating spline $\Gamma_{Q,\alpha}$, with $\alpha=0.5$, when a boundary 
point is perturbed, while part (b) shows the case when an interior point is perturbed. Local stability of the 
representation is apparent.\label{fig:shape_s} }
\end{figure}

Given a point cloud in two or three dimensions, many methods
of computational geometry are available to represent its boundary \cite{linsen2001point,pauly2003shape,pauly2004uncertainty,pauly2006point,berger2017survey}.
The structural simplicity of these representations allows efficient
handling for geometric deformations and topology changes. The results
of this paper are based on these ideas. In order to construct this
representation, we use the $\alpha$-shape algorithm that computes
the non-convex hull of a set of points in the plane. This algorithm
produces a polygonal hull ($\alpha$-shape) that depends on the value
of a real parameter (called $\alpha$). To obtain a smooth boundary,
we interpolate on the points that define the $\alpha$-shape with
a third order B-spline. B-splines are used due to its local stability
properties on the control points, see Figure \ref{fig:shape_s} (a),
(b). Palafox \emph{et al.} \cite{palafox2014effective} introduced
this procedure to estimate the boundary of a smooth sound obstacle
from far field data. A detailed description of the algorithm is as follows: \\

\begin{algorithm}[H]
\textbf{Data:} Let us $Q=\{p_{i}\}_{t=1}^{m}$ a set of $m$ points within the domain
of the problem, $Q\subset G\subset\mathbb{{R}}^{2}$, and $\alpha$
a positive real number.\\

\textbf{Output:} A third order B-spline using the point cloud non-convex hull as control
points.\\

\textbf{Step 1.} Compute the $\alpha$-shape of $Q$, denoted by $S_{\alpha}(Q)$,
using the $\alpha$-shape algorithm \cite{edelsbrunner1983shape}.

\textbf{Step 2.} If $S_{\alpha}(Q)$ is a valid polygon (i.e non-empty, non-self-intersected and non-disconnected): \\
Compute $\Gamma_{Q,\alpha}$ interpolating the points that define
the polygon $S_{\alpha}(Q)$ with a third-order B-spline.

\caption{Spline interpolation of the non-convex hull of a point cloud in 2D\label{alg:spline}}
\end{algorithm}

\begin{figure}[H]
\begin{centering}
\includegraphics[scale=0.25]{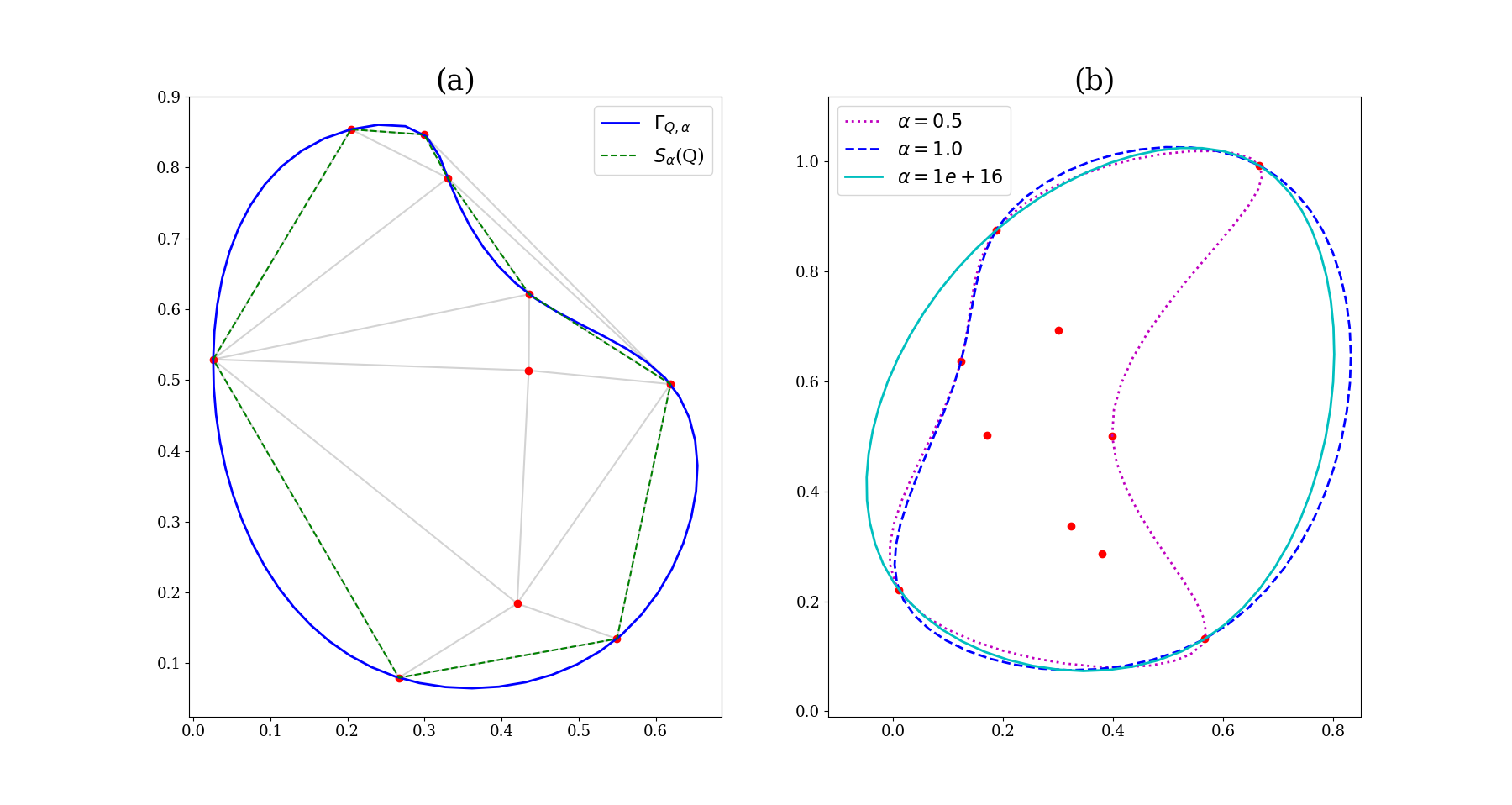}
\par\end{centering}
\caption{Part (a) shows the $\alpha$-shape (green), the interpolating spline (blue), $\Gamma_{Q,\alpha}$, with $\alpha=0.4$, and the correspondent Delaunay triangulation (grey). In part (b) shows $\Gamma_{Q,\alpha}$ for 
different values of $\alpha$.\label{fig:delanuay_shape}}
\end{figure}

Of note, the $\alpha$-shape algorithm uses the Delaunay triangulation with vertices Q to compute $S_{\alpha}(Q)$ (see Figure \ref{fig:delanuay_shape} (a)). 

The algorithm defined above accounts for a well defined mapping
from $Q$ and parameter $\alpha$ to a boundary $\Gamma_{Q,\alpha}$.
Indeed, the polygon $S_{\alpha}(Q)$, and thus the smooth boundary
$\Gamma_{Q,\alpha}$, is uniquely determined by $Q$. The positive parameter
$\alpha$ determines how non-convex the smooth boundary $\Gamma_{Q,\alpha}$
is (see Figure \ref{fig:delanuay_shape} (b)). 

\section{Bayesian Formulation of the Inverse Problem}\label{sec:bayes}

In the inverse scattering problem, the task is to estimate the shape of the scatterer
$D$ and the refractive index $\mathbf{b}$ given some available observations,
e.g., near field data observed at some points of the domain $G$.
The scattered field is obtained from several incoming plane waves
with different directions $d$. We shall formulate the inverse problem in a finite
dimensional space using the parametric representation of the obstacle
described above. Thereby, the parameter space to be estimated is $\theta=(Q,\alpha,{\bf b})$.

We shall pose the inverse scattering problem using the Bayesian approach.
This formulation allows us restating the inverse problem as ``...\emph{a well-posed
extension in a larger space of probability distributions}..." \cite{kaipio2006statistical}.
The solution of the inverse problem is the posterior probability distribution.\\

Let us consider the following additive noise pointwise observation
model: 
\begin{equation}
\mathbf{d}_{j}:=u(x_{j})+\eta_{j},\quad j=1,\ldots,M.,\label{eq:noise_point}
\end{equation}
where $\{x_{j}\}_{j=1}^{M}$, is the set of points at which the scattered
field is observed, $\eta_{x_{j}}\sim N(0,\sigma^{2})$ is additive
noise, and $\mathbf{d}_{j}$ noisy observations. Concatenating all
the observation, one can rewrite (\ref{eq:noise_point}) as 

\begin{equation}
\mathbf{d}:=F(\theta)+\eta\label{eq:model}
\end{equation}
with $F:=[u(x_{1}),\ldots,u(x_{M})]$ denoting the mapping from the
shape of $D$ (defined by Q and $\alpha$) and the refractive index
$\mathbf{b}$ to the noise free observables, $\eta$ being normally
distributed as $\mathcal{N}(0,\sigma^{2}I)$, and ${\bf d}=[\mathbf{d}_{1},\ldots,\mathbf{d}_{M}]^{T}$.

The Bayesian solution to our inverse problem is \cite{kaipio2006statistical}

\begin{equation}
\pi(\theta|\mathbf{d})=\frac{L(\mathbf{d}|\theta)\pi(\theta)}{\pi(\mathbf{d})}.\label{eq:posterior}
\end{equation}
where the likelihood is given by the noise distribution (\ref{eq:model}),
\begin{equation}
L(\mathbf{d}|\theta)=\left(2\pi\sigma\right)^{-M/2}\exp\left\{ -\frac{1}{2}\sum_{x_{j}}\left|\left|\frac{\mathbf{d}_{j}-u(x_{j})}{\sigma}\right|\right|^{2}\right\} \label{eq:likelihood}
\end{equation}
 assuming independent measurements. This function shows how measurements
affect knowledge about $\theta$.\\

The prior distribution expresses the information about $\theta$ independent
of the measurements $\mathbf{d}$. In our case, it is modeled as follows: 
\begin{enumerate}
\item The refractive index ${\bf b}$, the point cloud $Q$ and the shape
parameter $\alpha$ are taken as independent parameters; therefore
\[
\pi(\theta)=\pi(Q)\pi({\bf b})\pi(\alpha).
\]
\item Uniform distributions for $Q$ and $\alpha$ are imposed, as a default
selection. Consequently, the prior probabilities $\pi(Q)$, $\pi(\alpha)$
are constant on the domain $G$.
\item The choice of a gamma distribution for ${\bf b}$ is motivated by
the application of the inverse scattering problem in elasticity imaging
methods. The refractive index represents the Young modulus \cite{greenleaf2003selected}.
In this manner, the refractive index is assumed as a positive parameter.
\begin{equation}
{\bf b}\sim \text{Gamma}(\tilde{k},\tilde{\lambda}).\label{eq:Prior_b}
\end{equation}
\end{enumerate}

The resulting posterior distribution is explored using MCMC method
with Metropolis Hastings (MH) by generating a Markov chain with equilibrium
distribution (\ref{eq:posterior}). In the next Section we propose
a the probability transition kernel to implement the MH algorithm.

\subsection{An affine invariant probability transition kernel}\label{sec:kernel}

By construction, the performance of the Metropolis-Hastings (MH) algorithm
depends on a proposal probability distribution used to construct
the probability transition kernel. Often, proposal probability distributions
are gradient informed \cite{Girolami2011,bui2013computational}. However,
there is a class of derivative-free proposals, first introduced by
Christen and Fox \cite{christen2010}, for arbitrary continuous distributions
and requires no tuning, that gives rise to ``...\emph{an algorithm that is invariant to scale, and approximately invariant
to affine transformations of the state space}...''. This property implies that the performance of the MCMC method
is independent of the aspect ratio of anisotropic target distributions.

In Daza \emph{et al.}~\cite{daza2016solution}, a MCMC method with the affine invariant
property was used to approach the inverse scattering problem in an
inhomogeneous medium with a convex obstacle. In this paper we build
on previous results introducing an affine invariant proposal also
for the $\alpha$ parameter of the non-convex hull of a point cloud.
This proposal is defined from the circumradius of the Delaunay triangulation
of the point cloud at the iteration t. The circumradius is used for
two reasons: it contains information about the scale of the space
and the $\alpha$-shape is obtained with decision criteria that involve
the radius of the Delaunay triangulation (see details in \cite{edelsbrunner1983shape}).

The probability transition kernel for the parameter space $\theta^{\left(t\right)}=\left(Q^{(t)},b^{\left(t\right)},\alpha^{\left(t\right)}\right)$,
where $Q^{(t)}=\{p_{i}^{(t)}\}_{i=1}^{m}$ is given by

\[
K(\theta^{\left(t\right)},\theta^{\left(t+1\right)})=\sum_{i=1}^{4}w_{i}q_{i}(\theta^{\left(t\right)},\theta^{\left(t+1\right)}),
\]
and $w_{i}$ is the weight to choose the proposal $q_{i}$. Thus,
the new proposal is defined as follows:
\begin{enumerate}
\item \textbf{Point move:} Move randomly a single point $p_{k}^{(t)}\in Q^{(t)}=\{p_{i}^{(t)}\}_{i=1}^{m}$
\[
p_{k}^{(t+1)}=p_{k}^{(t)}+u_{-k}
\]
 where $u_{-k}\sim U(-\bar{d}_{-k},\bar{d}_{-k})$, and $\bar{d}_{-k}$
is the mean of the pairwise distances 
\[
|p_{j}^{(t)}-p_{l}^{(t)}|\hspace{2mm}\textsl{for}\hspace{2mm}j\neq l,\hspace{4mm}p_{j}^{(t)},p_{l}^{(t)}\in Q^{(t)}\setminus\{p_{k}^{(t)}\}.
\]
\item \textbf{Translate move:} Move randomly the point cloud
\[
p_{i}^{(t+1)}=p_{i}^{(t)}+u,\hspace{2mm}\textsl{for}\hspace{2mm}1<i<m
\]
 where $u\sim U(-\bar{d},\bar{d})$, and $\bar{d}$ is the mean of
the pairwise distances 
\[
|p_{j}^{(t)}-p_{l}^{(t)}|,\hspace{2mm}\textsl{for}\hspace{2mm}j\neq l\hspace{2mm},p_{j}^{(t)},p_{l}^{(t)}\in Q^{(t)}.
\]
\item \textbf{b move}: We move the refractive index ${\bf b}$ using its
prior distribution (\ref{eq:Prior_b}) 
\[
{\bf b}\sim\pi({\bf b}).
\]
\item \textbf{$\alpha$ move}: The move for the $\alpha$-shape parameter
is defined as follows:
\end{enumerate}
\begin{equation}
\alpha^{(t+1)}=\frac{1}{2}\alpha^{(t)}+\frac{1}{2}\mathcal{U}(r_{min},r_{max})\label{eq:prop_alpha-1}
\end{equation}
where $r_{min}$, $r_{max}$ are the minimum and maximum circumradius
of the Delaunay triangulation, respectively. The $\alpha$ value is
taken in $[r_{min},r_{max}]$, because for values of $\alpha\geq r_{max}$
the $\alpha$-shape of $Q$ is the convex hull and for values of $\alpha\leq r_{min}$
the $\alpha$-shape is the set of points $Q$.\\
\\
The instrumental proposals for the obstacle boundary (1-2) are symmetric
because when computing the pairwise distances without the point $p_{k}^{(t)}$
in the point move, we produce a symmetric proposal on the MCMC method. Similarly,
when the translation move is chosen, the entire cloud is moved with
the same direction $u$, and therefore the proposal is also symmetric
\cite{palafox2017point}. A detailed description of the algorithm
is as follows:\\

\begin{algorithm}[H]
\textbf{Data:} $\theta^{0}=(Q^{0},{\bf b}^{(0)},\alpha^{0})$. $Q^{(0)}=\{p_{i}^{(0)}\}_{i=1}^{m}$ uniformly distributed on $G$.
Number of samples $t_{max}$\\

\textbf{Output:} A prescribed number of samples of the posterior distribution \eqref{eq:posterior} obtained 
with the affine invariant sampler.\\

\textbf{Step 1.} Take $t=0$. Compute the energy 
\[
\text{Energy}(\theta^{(t)})=-\log\hspace{1mm}L(\mathbf{d}|\theta^{(t)})-\log\hspace{1mm}\pi({\bf b}^{(t)}).
\]

\textbf{Step 2.} Update $\theta^{(t+1)}$ randomly choosing one of the moves: \label{line:rand1}
\begin{enumerate}
\item Move a point $p_{k}^{(t)}\in Q^{(t)}$,
\[
p_{k}^{(t+1)}=p_{k}^{(t)}+u_{-k}
\]
\item Move every point in $Q^{(t)}=\{p_{i}^{(t)}\}_{i=1}^{m}$,
\[
p_{i}^{(t+1)}=p_{i}^{(t)}+u
\]
\item Draw ${\bf b}^{(t+1)}\sim \text{Gamma}(\tilde{k},\tilde{\lambda})$
\item Move $\alpha$, \\
\[
\alpha^{(t+1)}=0.5\alpha^{(t)}+0.5U(r_{min},r_{max})
\]
 \label{line:alpha_lin} 
\end{enumerate}
\textbf{Step 3.} Compute 
\[
\rho=\text{Energy}(\theta^{(t)})-\text{Energy}(\theta^{(t+1)})
\]
 and accept the proposal with probability $e^{-\rho}$. \\

\textbf{Step 4.} $t=t+1$\\

\textbf{Step 5.} Go back to \textbf{Step 2} if $t<t_{max}$, else stop.

\caption{Point Cloud Metropolis-Hastings Random Walk\label{alg:MH2}}
\end{algorithm}

\section{Numerical results and discussion}\label{sec:examples}

In the examples presented below we recover a star-shaper scatterer defined by
\begin{equation}
x(t)=(1.5\sin(t),\cos(t)+0.65\cos(2t)-0.65),\hspace{2mm}0\leq t\leq2\pi.\label{eq:kite}
\end{equation}

Of note, scatterer \eqref{eq:kite} was introduced by Kirsch \cite{kirsch1988two} 
and is often used in the inverse scattering literature as a benchmark. Recovering
scatterer \eqref{eq:kite} is challenging
for two reasons: first because it is non-convex and secondly, because
it possesses detailed structures (the wings) that are small compared
with the scale of the whole figure. In the setting for our numerical examples 
we take fixed refractive index ${\bf b}=25$. Furthermore, we use the area of the 
scatterer as a quantity of interest to explore the robustness of our inference method.

Data is measured at grid points of the numerical domain. In
order to avoid the inverse crime we solve the direct scattering problem
with the method introduced by Vainikko \cite{vainikko2000fast} to
create synthetic data. We add Gaussian noise
with mean zero and standard deviation $\sigma=0.012$, i.e., the signal
to noise ratio (SNR) is 100. 

We propose a design on the distribution of the
incident field over the unit circle. This design involves incident
fields in two wavenumbers ($k=1$, $k=5$). We take eight incident wave
directions uniformly distributed

\begin{equation}
\label{eq:inc_dirs}
d_{i}=\left(\cos\left(\frac{2\pi}{8}i+\zeta\right),\sin\left(\frac{2\pi}{8}i+\zeta\right)\right),\hspace{4mm}i=1,...,8.
\end{equation}

In directions $d_{2i}$, $i=1,\ldots,4$, the incident field
is applied with wavenumber $k=5$ and in the other directions with
wavenumber $k=1$. Data is assumed independent, then the likelihood
is given by

\begin{eqnarray*}
L(\mathbf{d}|\theta) & = & \left(2\pi\sigma_{L}\right)^{-M/2}\exp\left\{ -\frac{1}{2}\sum_{x_{j}}\left|\left|\frac{\mathbf{d}_{j}^{L}-u_{L}(x_{j})}{\sigma_{L}}\right|\right|^{2}\right\} +\\
 &  & +\left(2\pi\sigma_{H}\right)^{-M/2}\exp\left\{ -\frac{1}{2}\sum_{x_{j}}\left|\left|\frac{\mathbf{d}_{j}^{H}-u_{H}(x_{j})}{\sigma_{H}}\right|\right|^{2}\right\} ,
\end{eqnarray*}
where $u_{L}$ and $u_{H}$ denote the scattered field for the wavenumbers
$k_{L}=1$ and $k_{H}=5$ respectively. 

In Example \ref{exa:example1} we take $\zeta=0$ for the configuration 
of the incident field directions \eqref{eq:inc_dirs}. The mesh of the numerical domain is  
$N=40$ and $h=0.02$. In Example \ref{exa:example2} 
we explore the robustness of our method by rotating the direction of the incident fields 
with respect to Example \ref{exa:example1}. Finally, in the third example we verify that
by refining the mesh in the solution of the direct problem, the results improve accordingly. 
All numerical experiments were performed in python. For the sake of reproducibility, all code 
is available in \href{https://github.com/mdazatorres/MCMC_afin_invariant_sampling_ISP}{Github}.

\begin{example}
\label{exa:example1}
\end{example}

\begin{figure}
\begin{centering}
\includegraphics[scale=0.4]{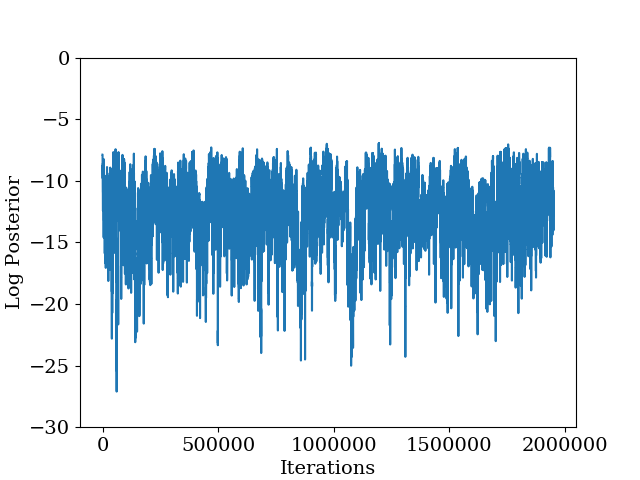}
\par\end{centering}
\caption{Trace plot of the logarithm of the posterior distribution with a burn-in of 50,000
iterations. It is apparent that the MCMC sampling is at equilibrium. \label{fig:energy_k15}}
\end{figure}

\begin{figure}
\begin{centering}
\includegraphics[scale=0.25]{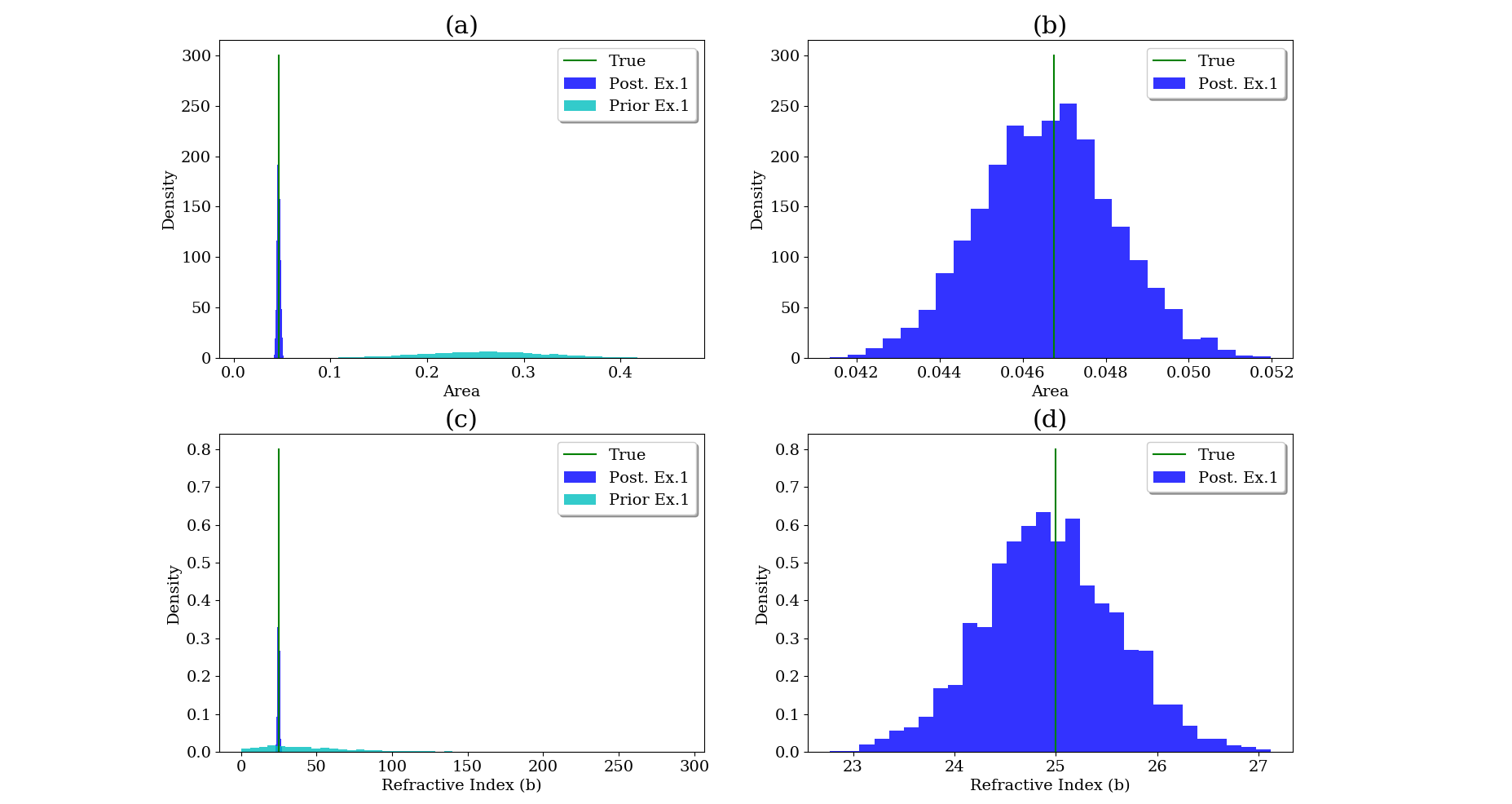}
\par\end{centering}
\caption{Parts (a) and (c) show the prior distribution in cyan, the posterior distribution in blue, and the true value
in green of both, the area of the scatterer and the refractive index for Example 1. Parts (b) and (d) show 
a close up along the horizontal axis near the true values. \label{fig:post_k1_post_k15}}
\end{figure}

We computed 2,000,000 samples of the posterior distribution. 
Figure \ref{fig:energy_k15} shows a trace plot of samples of the
logarithm of the posterior distribution at equilibrium, indicating
convergence in measure of the MCMC method. On the other hand, 
histograms of the prior and posterior distributions of both, the refractive 
index $\mathbf{b}$ and the area of the obstacle are shown in Figure 
\ref{fig:post_k1_post_k15}. True values are shown for comparison.

In order to show that the results obtained in Example \ref{exa:example1}
are independent of the directions of the incident field, in the next 
example we repeat the numerical experiment except that none of the incident 
directions coincide with an axis of symmetry of the obstacle. 

\begin{example}
\label{exa:example2}
\end{example}

\begin{figure}
\begin{centering}
\includegraphics[scale=0.25]{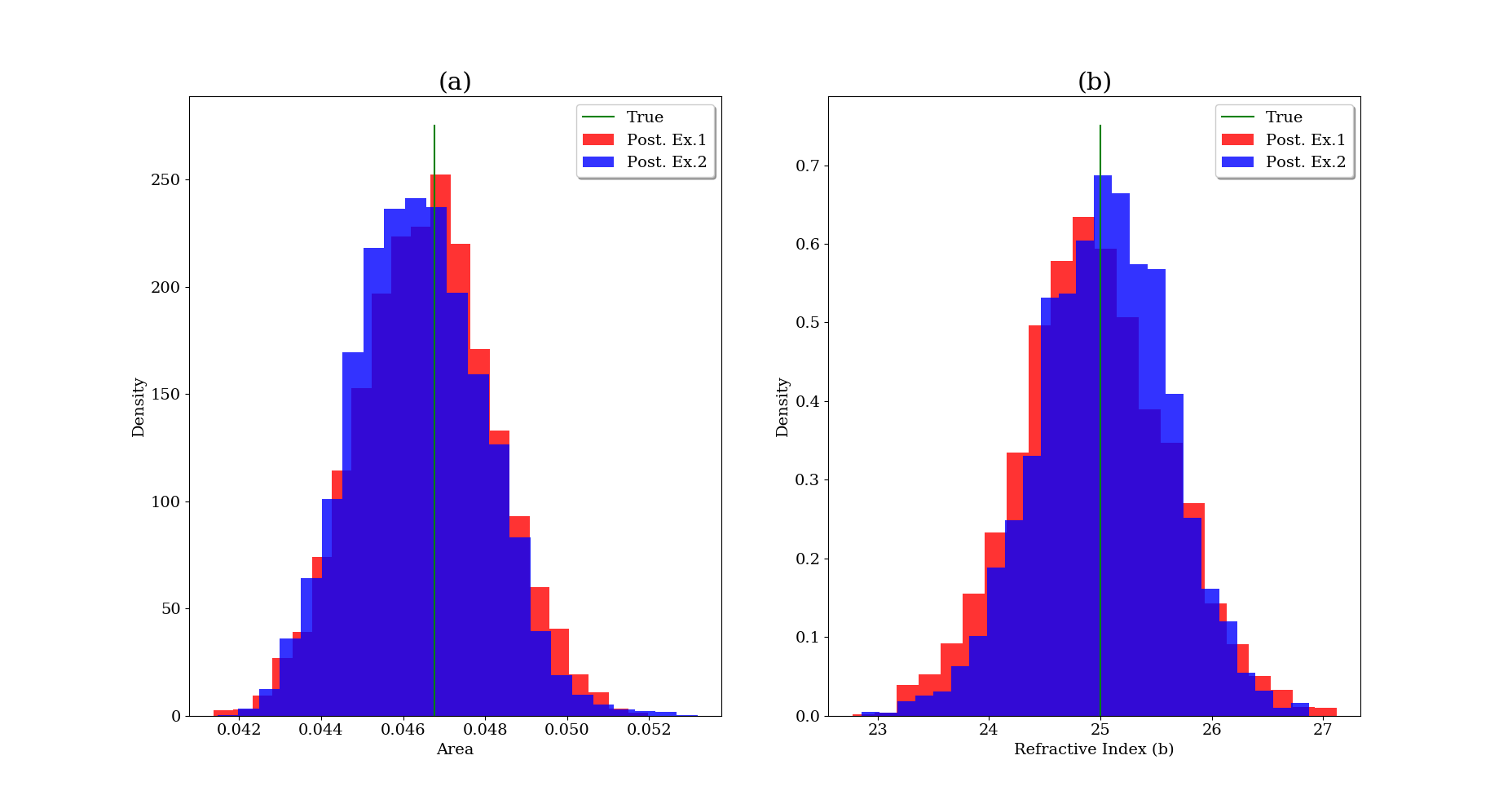}
\par\end{centering}
\caption{For examples \ref{exa:example1} and \ref{exa:example2}; posterior
distributions of the scatterer area and (b) the refractive index.
\label{fig:post_k15_post_k15r}}
\end{figure}

We consider the same setting of Example 1, except that the directions in
the incident field are rotated letting $\zeta=\frac{\pi}{6}$. 
Figure \ref{fig:post_k15_post_k15r} compares histograms of both,
the refractive index $\mathbf{b}$ and the area of the obstacle
obtained in Examples \ref{exa:example1} and \ref{exa:example2}.

Finally, Example \ref{exa:example3} illustrates the effect
of using a finer mesh in the numerical solution of the direct problem.
As expected, the variance of both, the refractive index $\mathbf{b}$
and the quantity of interest decreases if the mesh is reduced.

\begin{example}
\label{exa:example3}
\end{example}

\begin{figure}
\begin{centering}
\includegraphics[scale=0.3]{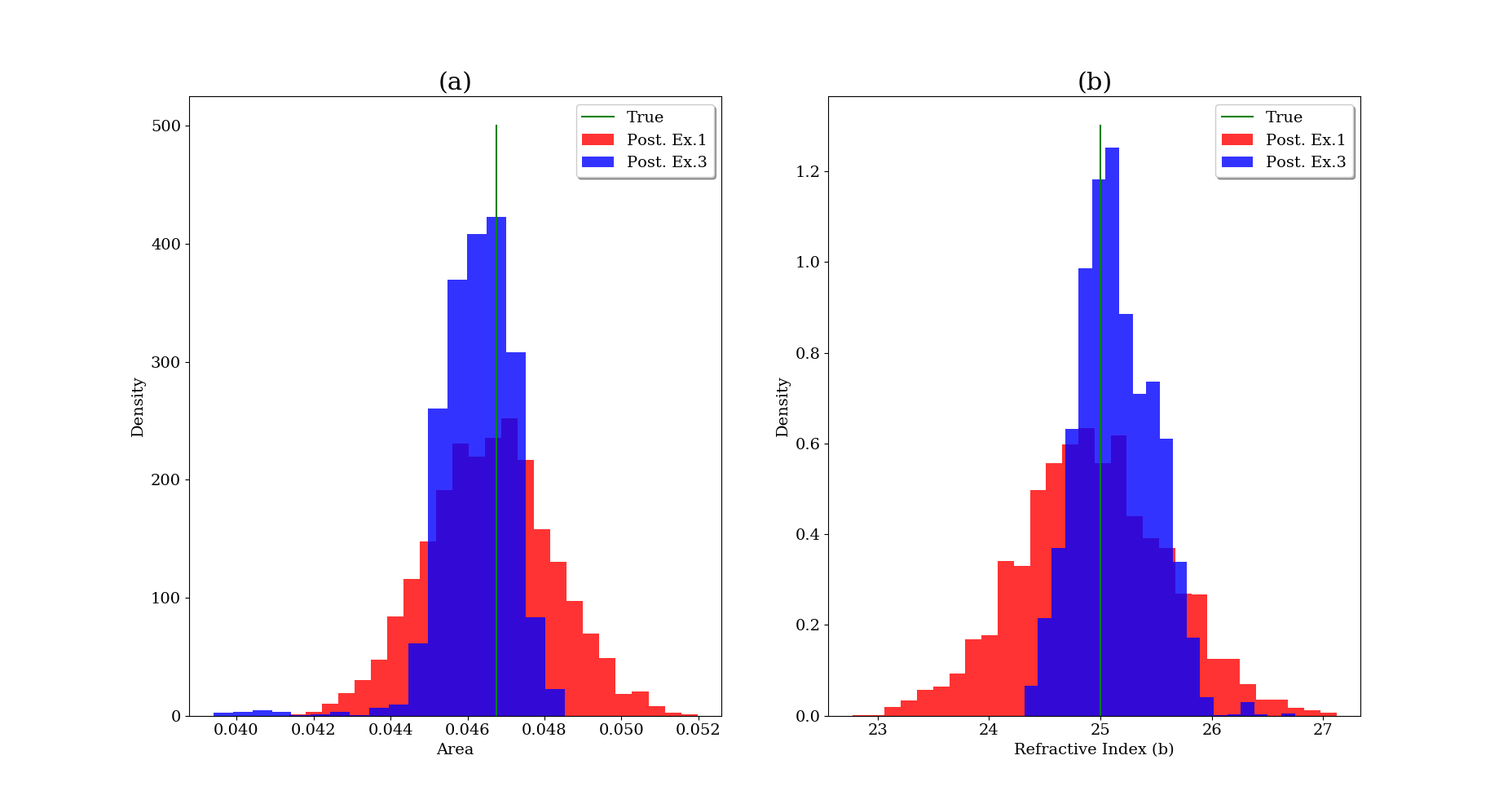}
\par\end{centering}
\caption{Example \ref{exa:example3}; posterior distributions of the scatterer area
and (b) the refractive index. \label{fig:post_k15_post_k15g}}
\end{figure}

\begin{figure}[H]
\centering
\begin{subfigure}{5cm}
\centering{}\includegraphics[scale=0.25]{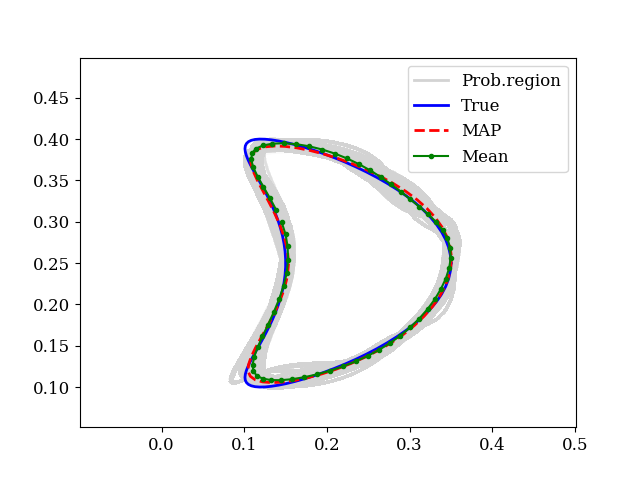}
\caption{Example \ref{exa:example1}}
\end{subfigure}%
\begin{subfigure}{5cm}
\centering\includegraphics[scale=0.25]{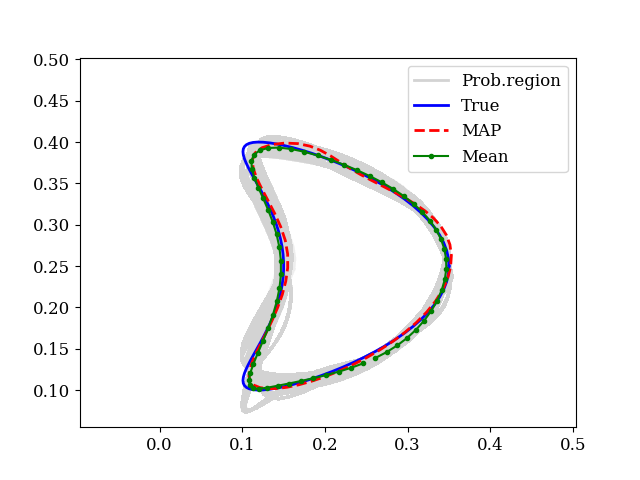}
\caption{Example \ref{exa:example2}}
\end{subfigure}\\
\begin{subfigure}{8cm}
\centering\includegraphics[scale=0.5]{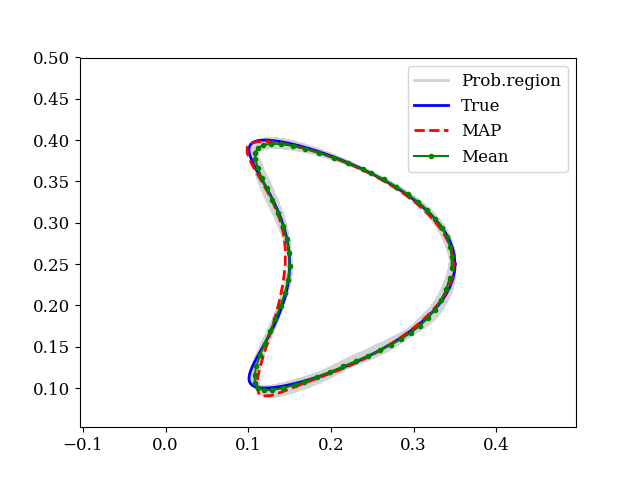}
\caption{Example \ref{exa:example3}}
\end{subfigure}
\caption{Numerical evidence in Examples \ref{exa:example1} and  \ref{exa:example2} suggests that the posterior 
distribution is insensitive with respect to a rotation of the direction of the incident fields. On the other hand, 
Example \ref{exa:example3} indicates there is a trade off between precision and 
computational cost as expected.\label{fig:map_cm_truth}}
\end{figure}

In this example, the configuration is taken as in Example \ref{exa:example1},
and we refined the mesh taking $N=80$, $h=0.01$. We computed 1,000,000
samples of the posterior distribution. Histograms of the posterior of both,
the refractive index $\mathbf{b}$ and the area of the obstacle are
shown in Figure \ref{fig:post_k15_post_k15g} together with the
respective results of Example \ref{exa:example1}. As before, the true values
are shown for comparison. Figure \ref{fig:map_cm_truth} summarizes the retrieval 
of the scatterer in all examples above. Numerical evidence 
suggest that the posterior distribution is insensitive with respect to a rotation 
of the direction of the incident fields. On the other hand, there is a trade off between 
precision and computational cost as expected. Table \ref{tab:estimators} compares
the true value, maximum a posteriori and conditional mean of our quantity of interest.
All estimators are correct with roughly two digits. Of note, the area of the true scatterer
was computed with the same resolution of the point cloud interpolation.

We may discuss the results of this manuscript as follows: we have introduced a method,
based on standard ideas of computational geometry, capable of simultaneously recovering
both, the support and constant refractive index of a star-shaped scatterer. In Figures
\ref{fig:energy_k15} through \ref{fig:map_cm_truth} and Table \ref{tab:estimators}, we offer numerical
evidence of the robustness of our method using a star-shaped scatterer commonly used as 
benchmark in the literature of inverse scattering problems under a particular data design. 

On the other hand, although our method does accomplish the task of reliably solving the inverse 
problem, it is computationally expensive. We remark that our method is derivative free. We believe that further 
research on affine invariant proposals is justified i.e., a natural extension is to explore the performance of 
derivative informed proposals of the probability transition kernel that fulfill the affine invariant property.
In general, we believe that it might be possible to further extend the idea of axiomatically imposing geometric
properties on move proposals for the probability transition kernel, such as moves that commute with
affine transformations of space, and then introducing examples of moves that fulfill those conditions.

\begin{table}[h]
\begin{center}
\begin{tabular}{|c|c|c|c|}
\cline{1-3}
\multicolumn{3}{ |c| }{Area of the scatterer}&\multicolumn{1}{ |c }{}\\
\cline{1-3}
True value&Conditional Mean&Maximum a posteriori&\multicolumn{1}{ |c }{}\\
\cline{1-4}
$4.68\times10^{-2}$&$4.66\times10^{-2}$&$4.60\times10^{-2}$&Example 1\\
\cline{1-4}
\multicolumn{1}{ c| }{}&$4.64\times10^{-2}$&$4.53\times10^{-2}$&Example 2\\
\cline{2-4}
\multicolumn{1}{ c| }{}&$4.63\times10^{-2}$&$4.51\times10^{-2}$&Example 3\\
\cline{2-4}
\end{tabular}
\end{center}
\caption{All estimators are correct with roughly two digits. Of note, the area of the true scatterer
was computed with the same resolution of the point cloud interpolation.\label{tab:estimators}}
\end{table}

\section{Conclusions}\label{sec:conclusions}

In this paper we have demonstrated a reliable method to approximate the solution
of the inverse scattering problem in the case of penetrable obstacles. Of note, we can 
recover simultaneously the support and constant refractive index of the scatterer. 
Our results are based on standard ideas of computational geometry, which  is a very active 
area of computer science with numerous applications. We consider that our findings are amenable 
to numerous extensions. 

We believe that future work should be directed towards improving the efficiency of the method
in terms of computational cost. Three possible research avenues, not necessarily in order of importance, 
are as follows. First, explore the performance of derivative informed proposals of the probability transition 
kernel that fulfill the affine invariant property. Second, defining a delayed acceptance method on the frecuency 
(lower wavenumber and higher wavenumber) spaces with two stages; the lower-model with a coarse mesh 
and the higher-model with a fine mesh. Finally, another approach might be, changing the number of points 
in the cloud in each level (delayed acceptance), for instance adding points in the cloud, taking more points 
for the interpolation of the point cloud.

\section*{Acknowledgements }

The authors would like to thank Colin Fox, Robert Scheichl, Tony Shardlow
and Antonio Capella for useful discussions and insight. J. A. I. R.
acknowledges the financial support of the Spanish Ministry of Economy
and Competitiveness under the project MTM2015-64865-P and the research
group MOMAT (Ref.910480) supported by \emph{Universidad Complutense de Madrid}.

\bibliographystyle{unsrt}
\bibliography{inverse_scattering}

\begin{thebibliography}{10}

\bibitem{Stuart2010}
A.~M. Stuart.
\newblock Inverse problems: A bayesian perspective.
\newblock {\em Acta Numerica}, 19:451--559, 2010.

\bibitem{christen2017posterior}
J~Andr{\'e}s Christen.
\newblock Posterior distribution existence and error control in banach spaces.
\newblock {\em arXiv preprint arXiv:1712.03299}, 2017.

\bibitem{bui2013computational}
Tan Bui-Thanh, Omar Ghattas, James Martin, and Georg Stadler.
\newblock A computational framework for infinite-dimensional bayesian inverse
  problems part i: The linearized case, with application to global seismic
  inversion.
\newblock {\em SIAM Journal on Scientific Computing}, 35(6):A2494--A2523, 2013.

\bibitem{Girolami2011}
Mark Girolami and Ben Calderhead.
\newblock Riemann manifold langevin and hamiltonian monte carlo methods.
\newblock {\em Journal of the Royal Statistical Society: Series B},
  73(2):123--214, 2011.

\bibitem{edelsbrunner1983shape}
Herbert Edelsbrunner, David Kirkpatrick, and Raimund Seidel.
\newblock On the shape of a set of points in the plane.
\newblock {\em IEEE Transactions on information theory}, 29(4):551--559, 1983.

\bibitem{birman1977estimates}
M~Sh Birman and Mikhail~Zakharovich Solomyak.
\newblock Estimates of singular numbers of integral operators.
\newblock {\em Russian Mathematical Surveys}, 32(1):15, 1977.

\bibitem{colton2012inverse}
David Colton and Rainer Kress.
\newblock {\em Inverse acoustic and electromagnetic scattering theory},
  volume~93.
\newblock Springer Science \& Business Media, 2012.

\bibitem{christen2010}
J~Andr{\'e}s Christen and Colin Fox.
\newblock A general purpose sampling algorithm for continuous distributions
  (the t-walk).
\newblock {\em Bayesian Analysis}, 5(2):263--281, 2010.

\bibitem{aguilar2002high}
Juan~C Aguilar and Yu~Chen.
\newblock High-order corrected trapezoidal quadrature rules for functions with
  a logarithmic singularity in 2-d.
\newblock {\em Computers \& Mathematics with Applications}, 44(8):1031--1039,
  2002.

\bibitem{vainikko2000fast}
Gennadi Vainikko.
\newblock Fast solvers of the lippmann-schwinger equation.
\newblock In {\em Direct and inverse problems of mathematical physics}, pages
  423--440. Springer, 2000.

\bibitem{linsen2001point}
Lars Linsen.
\newblock {\em Point cloud representation}.
\newblock Univ., Fak. f{\"u}r Informatik, Bibliothek, 2001.

\bibitem{pauly2003shape}
Mark Pauly, Richard Keiser, Leif~P Kobbelt, and Markus Gross.
\newblock Shape modeling with point-sampled geometry.
\newblock In {\em ACM Transactions on Graphics (TOG)}, pages 641--650. ACM,
  2003.

\bibitem{pauly2004uncertainty}
Mark Pauly, Niloy~J Mitra, and Leonidas~J Guibas.
\newblock Uncertainty and variability in point cloud surface data.
\newblock {\em SPBG}, 4:77--84, 2004.

\bibitem{pauly2006point}
Mark Pauly, Leif~P Kobbelt, and Markus Gross.
\newblock Point-based multiscale surface representation.
\newblock {\em ACM Transactions on Graphics (TOG)}, 25(2):177--193, 2006.

\bibitem{berger2017survey}
Matthew Berger, Andrea Tagliasacchi, Lee~M Seversky, Pierre Alliez, Gael
  Guennebaud, Joshua~A Levine, Andrei Sharf, and Claudio~T Silva.
\newblock A survey of surface reconstruction from point clouds.
\newblock In {\em Computer Graphics Forum}, pages 301--329. Wiley Online
  Library, 2017.

\bibitem{palafox2014effective}
Abel Palafox, Marcos~A Capistr{\'a}n, and J~Andr{\'e}s Christen.
\newblock Effective parameter dimension via bayesian model selection in the
  inverse acoustic scattering problem.
\newblock {\em Mathematical Problems in Engineering}, 2014, 2014.

\bibitem{kaipio2006statistical}
Jari Kaipio and Erkki Somersalo.
\newblock {\em Statistical and computational inverse problems}, volume 160.
\newblock Springer Science \& Business Media, 2006.

\bibitem{greenleaf2003selected}
James~F Greenleaf, Mostafa Fatemi, and Michael Insana.
\newblock Selected methods for imaging elastic properties of biological
  tissues.
\newblock {\em Annual review of biomedical engineering}, 5(1):57--78, 2003.

\bibitem{daza2016solution}
Maria~L Daza, Marcos~A Capistr{\'a}n, J~Andr{\'e}s Christen, and Lil{\'\i}
  Guadarrama.
\newblock Solution of the inverse scattering problem from inhomogeneous media
  using affine invariant sampling.
\newblock {\em Mathematical Methods in the Applied Sciences}, 2016.

\bibitem{palafox2017point}
Abel Palafox, Marcos Capistr{\'a}n, and J~Andr{\'e}s Christen.
\newblock Point cloud-based scatterer approximation and affine invariant
  sampling in the inverse scattering problem.
\newblock {\em Mathematical Methods in the Applied Sciences}, 40(9):3393--3403,
  2017.

\bibitem{kirsch1988two}
Andreas Kirsch and R~Kress.
\newblock Two methods for solving the inverse acoustic scattering problem.
\newblock {\em Inverse problems}, 4(3):749, 1988.

\end{thebibliography}

\end{document}